\documentclass[a4paper,12pt]{amsart}

\usepackage[
  top=2.5cm,       
  bottom=2.5cm,    
  left=3cm,       
  right=3cm       
]{geometry}

\usepackage{xcolor}
\usepackage{amsmath}
\usepackage{amssymb}
\usepackage{amsfonts}
\usepackage[T1]{fontenc}
\usepackage{libertine}
\usepackage[libertine]{newtxmath}
\usepackage{hyperref}
\usepackage{graphicx}
\usepackage{float}
\usepackage{amsthm}
\usepackage[nowidering]{yhmath}
\usepackage{enumerate}
\usepackage{xcolor}

\newtheorem{thm}{Theorem}
\newtheorem{prop}{Proposition}
\newtheorem{lemma}{Lemma}

\newtheorem{rem}{Remark}
\newtheorem{claim}{Claim}
\newtheorem{definition}{Definition}

\def\card{\textrm{card}}
\def\Tilde{\widetilde}
\def\diam{\operatorname{diam}}
\def\dist{\operatorname{dist}}

\title{The dimensions of curly and dipolar Besicovitch-Kakeya sets}
\author{Antonio C\'ordoba and Chuhan Sun}

\begin{document}

\begin{abstract}
    In this paper, we establish sharp Hausdorff dimension and box-counting dimension bounds for sets containing meridian circles and for sets whose difference sets contain a sphere.
    
    Fix two antipodal points $N$, $S$ on $S^{n-1}$, $n\ge3$, and consider the family of great circles passing through them. Assuming the Besicovitch–Kakeya overlapping conjecture in $\mathbb R^{n-1}$, we prove that every bounded set $K\subset\mathbb R^n$ containing a translated and dilated copy of each circle in this family satisfies $\dim_H K\ge n-1$. In particular, the same conclusion holds for bounded sets containing a circle in every orientation. This bound is sharp, as shown by the unit sphere. We also prove unconditionally that every bounded set $K\subset\mathbb R^n$, $n\ge2$, satisfying $S^{n-1}\subset K-K$ has lower box counting dimension greater than or equal to $n(n-1)/(2n-1)$, and the bound is also sharp.
\end{abstract}

\maketitle

\section{Besicovitch-Kakeya overlapping conjecture} In Euclidean space $\mathbb{R}^n$, let $\mathcal{B}_N^r$, $N\geq1$, be the collection of $n-$dimensional tubes of radius $r>0$, height $rN$, and arbitrary directions. 

Suppose that $\{R_1,R_2,\ldots,R_M\}\subset\mathcal{B}_N^r$ are tubes such that the angle between two of them is always greater than $\frac{1}{N}$. (In particular, this assumption implies that $M\lesssim N^{n-1}$). The following is an outstanding conjecture:

\textbf{Besicovitch-Kakeya overlapping conjecture}: Under the hypothesis above, we have the inequality:
\[
\int_{\mathbb{R}^n}\Big|\sum_{j=1}^M\chi_{R_j}\Big|^{\frac{n}{n-1}}\leq C(n)[\log(1+N)]^{a(n)}\sum_{j=1}^M|R_j|,
\]
where $C(n)$, $a(n)$ are finite constants depending only upon the dimension $n$.
\begin{rem}
         For many applications, it suffices to have the weaker version where $C(n)[\log (1+N)]^{a(n)}$ is replaced by $C_\varepsilon N^\varepsilon$ for every $\varepsilon>0$.
\end{rem}
\begin{rem}
    The proof of the two dimensional case of this conjecture is well-known, see \cite{A.Cordoba},\cite{A. C}. Recently, there is an OpenAI announcement of a three dimensional proof which is closely related to the work of Wang and Zahl \cite{Wang-Zahl}.
\end{rem}

\begin{rem}
     Recall that given a bounded set $E\subset\mathbb{R}^n$ and given $\alpha\ge0$ and $\delta>0$, we have the definition
     \begin{equation}\label{H}
     \mathcal{H}^\alpha_\delta(E):=\inf\Big\{\sum_{k}\left[\diam(Q_k)\right]^\alpha\Big\},
     \end{equation}
     where the infimum is taken over all countable coverings of E by cubes $Q_k$ of diameter less than $\delta$. Define the Hausdorff measure as
     \[
     \mathcal{H}^\alpha(E):=\lim_{\delta\rightarrow0}\mathcal{H}^\alpha_\delta(E).
     \]
     The limit above always exists, although it can be infinite, and has the following properties:
     \begin{enumerate}[(1)]
         \item 
         $\mathcal{H}^{\alpha}$ is countable additive, invariant under translations and rotations, and satisfies
         \[
         \mathcal{H}^\alpha(\lambda E)=\lambda^\alpha\mathcal{H}^\alpha(E), \textrm{ for }\lambda>0.
         \]
         \item
         $\mathcal{H}^0(E)$ counts the number of points inside E, while
         \[
         \mathcal{H}^n(E)=c_n|E|,
         \]
         where $c_n=\frac{\pi^{n/2}}{2^n(n/2)!}$, and $|E|$ denotes Lebesgue measure of $E$ in $\mathbb{R}^n$.
         \item 
         $\mathcal{H}^\alpha (E)<\infty$ implies that $\mathcal{H}^\beta(E)=0$ if $\beta>\alpha$. On the other hand, if $\beta<\alpha$ and $\mathcal{H}^\alpha(E)>0$, then $\mathcal{H}^\beta(E)=\infty$.
     \end{enumerate}
     Therefore, given a bounded set $E\subset \mathbb{R}^n$, there is a unique $\alpha$ such that $\mathcal{H}^\beta(E)=0$ if $\beta>\alpha$ and $\mathcal{H}^\beta(E)=\infty$ if $\beta<\alpha$.
     \begin{definition}
         The unique $\alpha$ is called the Hausdorff dimension of $E$. That is 
         \[
         \dim_H(E)=\sup\big\{\beta:\mathcal{H}^\beta(E)=\infty\big\}=\inf\big\{\alpha:\mathcal{H}^\alpha(E)=0\big\}.
         \]
     \end{definition}
     An important observation is that we can substitute cubes by balls or dyadic cubes belonging to a mesh $\Gamma_k$ of width $2^{-k}$, where $k$ is a positive integer, without changing the value of $\dim_H (E)$.
     
     There is another notion of dimension in $\mathbb{R}^n$, namely the box-counting or Minkowski dimension, $\dim_B(F)$, of a bounded set $F$, see \cite{Falconer}.
     We have several equivalent definitions, one of which is the following.
     
     Let $N_\varepsilon$ be the number of cubes in a mesh of side length $\varepsilon$ needed to cover $F$. Then
     \[
     \underline{\dim}_B(F)=\liminf_{\varepsilon\to0}\frac{\log N_\varepsilon}{|\log\varepsilon|}\quad\text{and}\quad\overline{\dim}_B(F)=\limsup_{\varepsilon\to0}\frac{\log N_\varepsilon}{|\log\varepsilon|}.
     \]
In the case that both limits above coincide, their common value is
the box-counting dimension $\dim_B(F)$.

If we denote by
$F_\varepsilon=\{x:\dist(x,F)\le\varepsilon\}$
the Minkowski sausage (following the terminology of
Mandelbrot~\cite{Mandelbrot}) of the set $F$, then $\dim_B(F)$
describes the way in which the measure $|F_\varepsilon|$
decreases, that is,
\[
\underline{\dim}_B(F)
=
n-\limsup_{\varepsilon\to0}
\left|\frac{\log|F_\varepsilon|}{\log\varepsilon}\right|
\quad\text{and}\quad
\overline{\dim}_B(F)
=
n-\liminf_{\varepsilon\to0}
\left|\frac{\log|F_\varepsilon|}{\log\varepsilon}\right|.
\]
Also we have the inequality:
\[
\underline{\dim}_B(F)\ge\dim_H(F).
\]
\end{rem}

\section{Curly Kakeya sets}
Let us consider a bounded set $K\subset\mathbb{R}^n$, $(n\geq3)$ containing a one-dimensional circle on each possible position in $\mathbb{R}^n$.
\begin{thm}\label{curly Kakeya 1}
    If the weak Besicovitch-Kakeya overlapping conjecture holds in $\mathbb{R}^{n-1}$, we have
    \[
    \dim_{H}(K)\geq n-1.
    \]
\end{thm}
\begin{rem}
    Observe that the unit sphere $S^{n-1}\subset\mathbb{R}^n$ is an example to conclude that the estimate above is sharp. The case $n=3$ has been proved in \cite{Antonio Cordoba} by a slightly different method, directly and without assuming the weak Besicovith-Kakeya overlapping conjecture. But, on the other hand, reference \cite{A. C} gives a proof that the whole Besicovitch-Kakeya overlapping conjecture holds in dimension $n=2$. 
\end{rem}

Theorem \ref{curly Kakeya 1} will be a corollary of the following stronger result:

Fix two antipodal points $N$, $S$ (the poles) on the unit sphere $S^{n-1}$, and let $\mathcal{E}$ denote the equator, defined as the intersection of $S^{n-1}$ with the hyperplane through the origin orthogonal to the line joining $N$ and $S$. Thus, $\mathcal{E}$ is the $S^{n-2}$ unit sphere. For each point $\xi\in\mathcal{E}$, let $\widetilde{C}_\xi$ be the great circle obtained by intersecting $S^{n-1}$ with the two-dimensional plane through $N$, $S$, and $\xi$. We may call $\{\widetilde{C}_\xi\}_{\xi\in\mathcal{E}}$ the family of meridian circles.

Suppose a bounded set $K\subset\mathbb{R}^n$ contains, for each $\xi\in\mathcal{E}$, a circle $C_{\xi,r(\xi)}$ of radius $r(\xi)>0$ obtained from $\widetilde{C}_\xi$ by dilation and translation. Then we have:
\begin{thm}\label{curly Kakeya 2}
Assuming that the weak Besicovitch-Kakeya overlapping holds in $\mathbb{R}^{n-1}$, we have:
    \[
    \dim_{H}(K)\geq n-1.
    \]
\end{thm}

\begin{proof}
    Let $\{Q_\sigma\}$ be a covering of the set $K$ by dyadic cubes of side length $\ell(Q_\sigma)\leq2^{-k_0}$, such that
    \[
    K\subset\bigcup_{\sigma}Q_\sigma=\bigcup_{k\ge k_0}\left(\bigcup_{\sigma}Q_\sigma^k\right),
    \]
    where the cubes $Q_\sigma^k$ belong to a dyadic mesh $\Gamma_k$ in $\mathbb{R}^n$, and $\ell(Q_\sigma^k)=2^{-k}$. Define
    \[
    E_k=\bigcup_{\sigma}Q_\sigma^k,\quad N_k=\card\{Q_\sigma^k\}.
    \]

    From the definition of Hausdorff dimension, to prove Theorem \ref{curly Kakeya 2}, it suffices to show that for any $\alpha<n-1$, there exists a positive $\lambda_\alpha>0$, such that
    \[
    \sum_{k\geq k_0}N_k 2^{-k\alpha}\geq \lambda_\alpha>0.
    \]
    In fact, it will be convenient to prove a stronger result, namely:
    \begin{claim}
        For every $\alpha<n-1$, there exsits $\lambda_\alpha>0$, $k\geq k_0$ such that
        \[
        N_k2^{-k\alpha}\geq\lambda_\alpha>0.
        \]
    \end{claim}
    For each $\xi\in\mathcal{E}$, let $r(\xi)$ denote the radius of the circles $C_\xi\subset K$ obtained from $\widetilde{C}_\xi$ by dilation and translation.
    
    Let 
    \[
    \mathcal{E}_m=\{\xi\in\mathcal{E}:2^{m-1}\leq r(\xi)<2^m\},\quad m\in\mathbb{Z},
    \]
    and observe that
    \[
    \mathcal{E}=\bigcup_m\mathcal{E}_m
    \]
    implying that for some $m\in\mathbb{Z}$, we must have 
    \[
    \mathcal{H}^{\ast(n-2)}(\mathcal{E}_m)>0,
    \]
    where $\mathcal{H}^{\ast(n-2)}$ denotes the exterior $\mathcal{H}^{n-2}-$Hausdorff measure.
    Taking into consideration that Hausdorff dimension is invariant under dialtions, there is no loss of generality if we assume that $\mathcal{H}^{\ast(n-2)}(\mathcal{E}_1)>0$. Consequently, in the following we shall assume that $K$ contains a circle of radius $r(\xi)$ between $1$ and $2$ for each direction $\xi\in\mathcal{E}_1$, and  $\mathcal{H}^{\ast(n-2)}(\mathcal{E}_1)>0$.

    Let $C_\xi$ be one of those circles, our assignment of a North-South line $N-S$ yields automatically two poles $N_\xi-S_\xi$ in $C_\xi$. Let us divide $C_\xi$ into eight equal octants $\{C_\xi^j\}_{j=1}^{8}$, and choose one of them, says $C_\xi^1$ containing the north-pole $N_\xi$ as its middle point.

    \begin{lemma}
        For every $\beta>0$, there exists $\lambda_\beta>0$ such that, for each $\xi\in\mathcal{E}_1$,
        \[
        \mathcal{H}^1(C_\xi^1\cap E_k)\geq\lambda_\beta 2^{-k\beta},
        \]
        for a certain $k\geq k_0$.
    \end{lemma}
    \begin{proof}
        Assume there exists a $\beta>0$, such that
        \[
        \mathcal{H}^1(C_\xi^1\cap E_k)<\lambda_\beta 2^{-k\beta},\quad \textrm{ for any } \lambda_\beta>0.
        \]
        Then we have
        \[
        \frac{\pi}{4}\leq\mathcal{H}^1(C_\xi^1)=\sum_{k\geq k_0}\mathcal{H}^1(C_\xi^1\cap E_k)<\lambda_\beta\sum_{k\geq k_0}2^{-k\beta}\leq\lambda_\beta\frac{2^{-k_0\beta}}{1-2^{-\beta}},
        \]
        which is false for $\lambda_\beta$ small enough. Therefore, for each $\xi\in\mathcal{E}_1$, we have
        \[
        \mathcal{H}^1(C_\xi^1\cap E_k)\geq\lambda_\beta 2^{-k\beta}, \quad\textrm{  for a certain } k\ge k_0.
        \]
    \end{proof}
    Next, let us define the sets
    \[
    \Omega_k=\{\xi\in\mathcal{E}_1: \mathcal{H}^1(C_\xi^1\cap E_k)\geq \lambda_\beta 2^{-k\beta}\},
    \]
    and observe that
    \[
    \mathcal{E}_1\subset\bigcup_{k\geq k_0}\Omega_k.
    \]
    Since for each $\xi\in\mathcal{E}_1$, there exists $k\ge k_0$ satisfing $\xi\in\Omega_k$, so we have
    \[
    \mathcal{H}^{\ast(n-2)}(\Omega_k)\geq \lambda_\beta 2^{-k\beta},
    \]
    for a convenient $\lambda_\beta>0$.
    From this estimate we can conclude that the set $\Omega_k$ contains a subset of directions $\mathcal{A}=\{\xi_{a}\}$ such that $\left|\xi_a-\xi_b\right|\geq 2^{-k}$ for different $\xi_a$, $\xi_b\in\mathcal{A}$ and $\card\{\mathcal{A}\}\simeq\lambda_\beta 2^{-k\beta}2^{k(n-2)}$.
    Let us denote by $D_\xi^k$ the union of all the cubes $\{Q_\sigma^k\}$ in $E_k$ such that $C_\xi^1\cap Q_\sigma^k\neq \varnothing$ and observe that
    \[
    \left|D_\xi^k\right|\gtrsim\lambda_\beta2^{-k\beta}2^{-k(n-1)}.
    \]
    Then we have
    \begin{equation}\label{ineq}
    \sum_{\xi_a\in\mathcal{A}}\left|D_{\xi_a}^k\right|\gtrsim\lambda_\beta2^{-k\beta}2^{-k(n-1)}\lambda_\beta2^{-k\beta}2^{k(n-2)}=\lambda_\beta^22^{-k}2^{-2k\beta},
    \end{equation}
    and H\"older's inequality yields
    \begin{equation}\label{star}
        \sum_{\xi_a\in\mathcal{A}}|D_{\xi_a}^k|\leq\left|E_k\right|^{\frac{1}{n-1}}\left[\int\left|\sum_{\xi_a\in\mathcal{A}}\chi_{D_{\xi_a}^k}\right|^{\frac{n-1}{n-2}}\right]^{\frac{n-2}{n-1}}.
    \end{equation}
    Next we project every $D_{\xi_a}^k$ onto a set $\Tilde{D}_{\xi_a}^k$ by a vertical projection so that each cube $Q_\sigma^k$ satisfying $Q_\sigma^k\cap C_{\xi_a}^1\neq\varnothing$ is translated vertically to the cube $\Tilde{Q}_\sigma^k$ centered at the hyperplane $x_n=0$ of $\mathbb{R}^n$. Let us call $\Tilde{D}_{\xi_a}^k$ the corresponding projected figure. 
    \begin{figure}[htbp]
    \centering
    \includegraphics[width=0.6\linewidth]{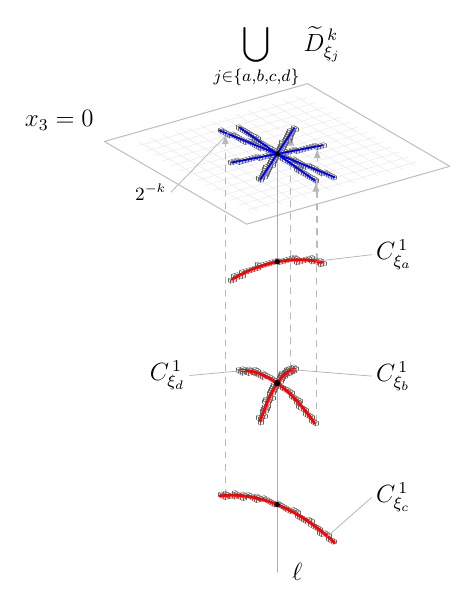}
    \caption{3-dimensional case}
    \label{figure}
    \end{figure}
    It happens that $\Tilde{D}_{\xi_a}$ is contained in a tube $R_{\xi_a}$ in $\mathbb{R}^n$ having vertical dimension $2^{-k}$ and directions $\xi_a$ in the $(n-1)$-horizontal hyperplane which are $2^{-k}$-separated as in the $(n-1)$-dimensional weak Besicovitch-Kakeya overlapping conjecture. Then an intergration in both the vertival and horizonal variables shows that the integral in the inequlity \eqref{star} is controlled by the corresponding integral
    \begin{equation}\label{star2}
       \int\left|\sum_{\xi_a\in\mathcal{A}}\chi_{R_{\xi_a}}\right|^{\frac{n-1}{n-2}}.
    \end{equation}
    To see this, let us denote by $Q_{\omega,j}$ the cubes in the covering $\{Q_\sigma^k\}$ having non-empty intersection with at least one of the class $C_{\xi_a}^1$ where $\omega$ is the index for the horizontal coordinates of $Q_{\omega,j}$ in the mesh $\Gamma_{k,j}$ and $j$ describes the vertical position also in $\Gamma_k$.
    The integral \eqref{star} restricted to the vertical band of $\Gamma_k$ describe by $\omega$ will be equal to a sum $2^{-kn}\sum_{\ell}\nu_\ell^{\frac{n-1}{n-2}}$ where the index $\ell$ describes cubes in the band intersecting $\nu_\ell$ arcs of the family $\{C_{\xi_a}^1\}$. On the other hand, we corresponding horizontal-vertical splitting of the integral \eqref{star2} yields $2^{-kn}\left(\sum_\ell \nu_\ell\right)^{\frac{n-1}{n-2}}$ for the same band, and we have the inequality
    \[
    \left(\sum_\ell \nu_\ell\right)^{\frac{n-1}{n-2}}\geq \sum_\ell \nu_\ell^{\frac{n-1}{n-2}}.
    \]
    Therefore inequalities \eqref{ineq} and \eqref{star} yields
    \[
    \lambda_\beta^22^{-k}2^{-2k\beta}\lesssim\left|E_k\right|^{\frac{1}{n-1}}\left[\int\left|\sum_{\xi_a\in\mathcal{A}}\chi_{R_{\xi_a}}\right|^{\frac{n-1}{n-2}}\right]^{\frac{n-2}{n-1}}.
    \]
    For any $\gamma>0$, as the weak Besicovitch-Kakeya overlapping estimate implies
    \[
    \left|E_k\right|\geq2^{-k\gamma}2^{-k}2^{-\mu k\beta},
    \]
    for some universal number $\mu$.
    That is
    \[
    N_k2^{-k\alpha}\geq2^{k[n-1-\alpha-\mu\beta-\gamma]}\geq\lambda>0,
    \]
    for $\alpha<n-1$ if $\beta$ and $\gamma$ are small enough.
\end{proof}

\section{dipolar Kakeya set}
 Let us consider a set $K\subset\mathbb{R}^n$, $n\ge 2$, such that $S^{n-1}\subset K-K$. That is, for each $\theta\in S^{n-1}$, there are two points $z(\theta), w(\theta)\in K$ such that $z(\theta)-w(\theta)=\theta$. The set $K$ contains both endpoints of a unit segment in every direction, instead of the entire segment. We call such a set $K$ a dipolar Besicovitch–Kakeya set.

\begin{thm}
    Let $K\subset\mathbb{R}^n$, $n\geq2$, be a bounded dipolar Besicovitch-Kakeya set. Then
    \[
    \underline{\dim}_{B}(K)\geq \frac{n(n-1)}{2n-1}.
    \]
    Moreover, the lower bound is sharp. 
\end{thm}
 
\begin{rem}
If the distance between the two points is allowed to vary with the direction, so that for each $\theta\in S^{n-1}$ there exist $z(\theta),w(\theta)\in K$ satisfying $z(\theta)-w(\theta)=r(\theta)\theta$ with $1\le r(\theta)\le2$, then $\underline{\dim}_B K\ge (n-1)/2$. See Proposition 5 in \cite{Antonio Cordoba} for the $n=2$ case, and the same covering argument works in all dimensions.
\end{rem}

\begin{proof}
    Let $E_k=\{Q_\sigma\}$ be a covering of $K$ by squares of a mesh $\Gamma_k$ in $\mathbb{R}^n$ with length $2^{-k}$, i.e. 
    \[
    \ell(Q_\sigma)=2^{-k}, \quad K\subset\bigcup Q_\sigma,
    \]
    and define $N_k=\textrm{card}(E_k)$.
    
    We fix a set of directions
    \[
    A_k=\{\theta_j\}_{j=1,2,\cdots,2^{k(n-1)}}\subset S^{n-1},
    \]
    which are at distance between $2^{-k}$ and $2^{-k+1}$ in the unit sphere.
    Define
    \[
    z(\theta_j)\in Q^z(\theta_j)\in E_k\quad\textrm{and}\quad w(\theta_j)\in Q^w(\theta_j)\in E_k,
    \]
    with
    \[
    z(\theta_j)-w(\theta_j)=\theta_j.
    \]
    Let
    \begin{align*}
        D(\theta_j)&=Q^z(\theta_j)\cup Q^w(\theta_j).\\
        \mathcal{F}&=\{\{Q^z(\theta_j), Q^w(\theta_j)\}: 1\leq j\leq 2^{k(n-1)}\}.\\
        \mathcal{R}(Q)&=\{Q_\sigma\in E_k:\{Q, Q_\sigma\}\in\mathcal{F}\}.\\
        M(Q)&=\textrm{card}\{\mathcal{R}(Q)\}=\textrm{the number of directions } \theta_j\in A_k, \textrm{ so that } D(\theta_j)\supset Q. 
    \end{align*}
    Next let us choose $Q_1$ such that $M(Q_1)$ is maximum. Having chosen $Q_1, Q_2,\cdots,Q_\tau$, we pick $Q_{\tau+1}$ with the same criterion, namely, $M(Q_{\tau+1})$ is maximum in $E_k\backslash\{Q_1,Q_2,\ldots,Q_\tau\}$. We continue this selection process until all the directions in $A_k$ are obtained within the sets
    \[
    \{Q_1,\cdots,Q_T\}\cup\bigcup_{j=1}^T\{Q_\sigma: Q_\sigma\in\mathcal{R}(Q_j)\}.
    \]
    We have
    \[
    \sum_{j=1}^TM(Q_j)\simeq2^{k(n-1)},
    \]
    and from this we conclude the existence of $\nu\in\{1, 2, \cdots, k(n-1)\}$ such that
    \begin{equation}\label{M nu 2}
    \frac{2^{k(n-1)}}{k(n-1)}\leq M_\nu~2^\nu\leq2^{k(n-1)},
    \end{equation}
    where
    \[
    M_\nu=\textrm{card}\{Q_j: 2^{\nu-1}\leq M(Q_j)<2^\nu\}.
    \]
    Let us denote by $\Tilde{Q}_1, \cdots,\Tilde{Q}_{M_\nu}$ those squares in the family $\{Q_1,\cdots,Q_T\}$ satisfying the estimate
    \[
    2^{\nu-1}\leq\textrm{card}(\mathcal{R}(\Tilde{Q}_j))<2^\nu.
    \]
    We have
    \begin{equation}\label{Est}
    \frac{2^{-k}}{(n-1)k}\lesssim M_\nu~2^\nu2^{-kn}\lesssim\sum_{j=1}^{M_\nu}\left(\sum_{Q\in\mathcal{R}(\Tilde{Q}_j)}|Q|\right).
    \end{equation}
    
     If $M_\nu\geq 2^{\frac{n(n-1)}{2n-1}k}$, then
    \[
    |E_k|=N_k2^{-kn}\geq M_\nu 2^{-kn}\geq 2^{\frac{n(n-1)}{2n-1}k}2^{-kn},
    \]
    which implies 
    \[
    N_k\geq 2^{\frac{n(n-1)}{2n-1}k},
    \]
    Then 
    \begin{equation}\label{case 1}
    \underline{\dim}_B(K)=\liminf_{k\rightarrow{\infty}}\frac{\log N_k}{|\log 2^{-k}|}\geq \liminf_{k\rightarrow{\infty}}\frac{\frac{n(n-1)}{2n-1}k}{k}=\frac{n(n-1)}{2n-1}.
    \end{equation}

     If $M_\nu <2^{\frac{n(n-1)}{2n-1}k}$, we estimate the sum in \eqref{Est} by H\"older inequality, 
    \[
    \sum_{j=1}^{M_\nu}\left(\sum_{Q\in\mathcal{R}(\Tilde{Q}_j)}|Q|\right)=\int\sum_{j=1}^{M_\nu}\sum_{Q\in\mathcal{R}(\Tilde{Q}_j)}\chi_Q\leq |E_k|^{\frac{n-1}{n}}\left(\int\left[\sum_{j=1}^{M_\nu}\sum_{Q\in\mathcal{R}(\Tilde{Q}_j)}\chi_{Q}\right]^{n}\right)^{\frac{1}{n}}.
    \]
    We can split the inner sum into $9^n$ parts, so that in each one of them the corresponding cubes $\Tilde{Q}_j$ are placed at distance greater than $8\times 2^{-k}$. Then Minkowski inequality will allow us to estimate the integral above under that assumption on the distance among $\Tilde{Q}_j$. Therefore we obtain
    \[
    \left(\int\left[\sum_{j=1}^{M_\nu}\sum_{Q\in\mathcal{R}(\Tilde{Q}_j)}\chi_{Q}\right]^n\right)^{\frac{1}{n}}\leq\sum_{\sigma=1}^{9^n}\left(\int\left[\sum_{j_\sigma=1}^{M_{\nu}^{\sigma}}\sum_{Q\in\mathcal{R}(\Tilde{Q}_{j_{\sigma}})}\chi_Q\right]^n\right)^{\frac{1}{n}},
    \]
    where $M_\nu^\sigma\leq M_\nu$, and for any $1\leq j_\sigma,j_{\sigma}'\leq M_{\nu}^\sigma$, $j_\sigma\neq j_{\sigma}'$, the distance between  $\Tilde{Q}_{j_\sigma}$ and $\Tilde{Q}_{j_\sigma'}$ is greater than $8\times 2^{-k}$.

    For each $\sigma$, we have
    \[
    \int\left[\sum_{j_{\sigma}=1}^{M_\nu^\sigma}\sum_{Q\in\mathcal{R}(\Tilde{Q}_{j_\sigma})}\chi_{Q}\right]^{n}\leq\sum_{j_{\sigma_1},j_{\sigma_2},\ldots,j_{\sigma_n}=1}^{M_\nu}\int\sum_{\substack{Q^1\in\mathcal{R}(\Tilde{Q}_{j_{\sigma_1}}) \\ Q^2\in \mathcal{R}(\Tilde{Q}_{j_{\sigma_2}})\\\vdots\\ Q^n\in\mathcal{R}(\Tilde{Q}_{j_{\sigma_n}})}}\chi_{Q^1}\chi_{Q^2}\ldots\chi_{Q^n},
    \]
    and observe that because the family $\mathcal{R}(Q)$ is disjoint, we have
    \[
    \left(\sum_{Q^1\in\mathcal{R}(Q)}\chi_{Q^1}\right)^\ell=\sum_{Q^1\in\mathcal{R}(Q)}\chi_{Q^1}.
    \]
    Let us write $X_m=\Big\{\{j_{\sigma_1},j_{\sigma_2},\dots,j_{\sigma_n}\}:\card\{j_{\sigma_1},j_{\sigma_2},\dots,j_{\sigma_n}\}=m\Big\}$, then
    \[
    \card \{X_m\}\leq M_\nu^m,
    \]
    and 
    \begin{equation}\label{exponent est}
    \sum_{\{j_{\sigma_1},j_{\sigma_2},\dots,j_{\sigma_n}\}\in X_m}\int\sum_{Q^1\in\mathcal{R}(\Tilde{Q}_{j_{\sigma_1}})}\cdots\sum_{Q^n\in\mathcal{R}(\Tilde{Q}_{j_{\sigma_n}})}\chi_{Q_1}\cdots\chi_{Q_n}\lesssim M_{\nu}^m2^{k(n-m)}2^{-kn}=(M_\nu2^{-k})^m\leq2^{\left[\frac{(n^2-3n+1)m}{2n-1}\right]k},
    \end{equation}
    because for $m$ different cubes in $\{\Tilde{Q}_{j_{\sigma_1}}, \Tilde{Q}_{j_{\sigma_2}}, \dots, \Tilde{Q}_{j_{\sigma_n}}\}$, $2\leq m\leq n$, the points at unit distance from the center of them lie on a sphere ${S}^{n-m}$, and there are at most $C_{n,m}2^{k(n-m)}$ cubes $Q$ satisfying $Q\in\bigcap_{i=1}^n\mathcal{R}(\Tilde{Q}_{j_{\sigma_i}})$.

    Furthermore, when $j_{\sigma_1}=j_{\sigma_2}=\ldots=j_{\sigma_n}=j$, using the condition \eqref{M nu 2} we have the estimate
    
    \begin{equation}\label{j}
    \sum_{j=1}^{M_\nu}\int\sum_{\substack{Q^1\in\mathcal{R}(\Tilde{Q}_{j}) \\ Q^2\in \mathcal{R}(\Tilde{Q}_{j})\\\vdots\\ Q^n\in\mathcal{R}(\Tilde{Q}_{j})}}\chi_{Q^1}\chi_{Q^2}\ldots\chi_{Q^n}
    \lesssim M_\nu 2^{\nu}2^{-kn}\lesssim 2^{-k}\leq 2^{\frac{(n^2-3n+1)n}{2n-1}k}.
    \end{equation}
    The last inequality above follows from the estimate
    \[
    \frac{(n^2-3n+1)n}{2n-1}\geq -\frac{2}{3}>-1, \quad \textrm{for } n\geq2.
    \]


    For $n\geq 3$, since $m\leq n$, the exponent in \eqref{exponent est} satisfies
    \[
    \frac{(n^2-3n+1)m}{2n-1}k\leq \frac{(n^2-3n+1)n}{2n-1}k.
    \]

    Therefore we obtain the estimate
    
    \[
    \frac{2^{-k}}{(n-1)k}\lesssim|E_k|^{\frac{n-1}{n}}2^{\frac{n^2-3n+1}{2n-1}k},
    \]
    and this implies
    \[
    |E_k|\gtrsim\frac{2^{-\frac{n^2}{2n-1}k}}{[(n-1)k]^{\frac{n}{n-1}}}.
    \]
    So we have
    \[
    \underline{\dim}_B(K)\geq n-\limsup_{k\rightarrow{\infty}}\left|\frac{\log|E_k|}{\log2^{-k}}\right|\geq n-\limsup_{k\rightarrow{\infty}}\frac{\frac{n^2}{2n-1}k+\frac{n}{n-1}\log [(n-1)k]}{k}=\frac{n(n-1)}{2n-1}.
    \]
    For the case $n=2$, we get 
    \[
    |E_k|\gtrsim\frac{2^{-\frac{4}{3}k}}{k^{2}},
    \]
    implying
    \[
    \underline{\dim}_B(K)\geq\frac{2}{3}.
    \]
    
To show the lower bound is sharp, we need the following proposition:
\begin{prop}\cite[Section~4.1]{Shmerkin-Yu}\label{example}
    For every integer $n\ge2$, there exists a compact set $K\subset\mathbb R^n$ satisfying $S^{n-1}\subset K-K$ such that 
\[
\underline{\dim}_{B}(K)\le\frac{n(n-1)}{2n-1}.
\]
\end{prop} 
For $n=2$, the detailed proof of Proposition \ref{example} is presented in \cite{Shmerkin-Yu}, Section 4.1. For the sake of completeness, we discuss in the following the general n-dimensional case.
\begin{proof}
We construct a sequence of finite point sets together with families of spherical caps. At each step, we
replace each cap by smaller caps clustered near the center of its supporting sphere, while preserving all of its directions up to sign. This clustering will allow us to obtain an efficient covering at a suitable sequence of scales.

Denote the unit sphere centered at $e \in \mathbb{R}^n$ by
\[
S(e,1)=\{x\in\mathbb R^n:|x-e|=1\}.
\]
And denote the closed ball of radius $r$ centered at $e\in\mathbb{R}^n$ by
\[
B(e,r)=\{x\in\mathbb{R}^n:|x-e|\leq r\}.
\]
For a set $F\subset\mathbb R^n$, denote its closed $\rho$-neighbourhood by
\[
F^{(\rho)}=\{x:\operatorname{dist}(x,F)\le\rho\}.
\]

 Let $E\subset S(e,1)$ be a spherical cap, and let $\delta>0$. Choose a finite maximal $\delta$-separated set
\[
\mathcal P(E,\delta)\subset E.
\]
The points in $\mathcal P(E,\delta)$ are pairwise at least $\delta$ apart, and maximality ensures that every point of $E$ lies within distance $\delta$ of some point in $\mathcal P(E,\delta)$.

For each $p\in\mathcal P(E,\delta)$, let
\[
E'_p=B(p,2\delta)\cap S(e,1).
\]
$E_p'$ is the part of the supporting sphere within distance
$2\delta$ of $p$.

Reflect $E_p'$ about the midpoint of the segment joining $e$ to $p$, and denote the reflected cap by
\[
E_p=e+p-E'_p
=\{e+p-z:z\in E'_p\}.
\]
The map $z\mapsto e+p-z$ reflects the unit sphere centered at $e$ to the
unit sphere centered at $p$, and reflects the ball centered at $p$ to the ball centered at $e$. Thus
\[
E_p=B(e,2\delta)\cap S(p,1).
\]
In particular, the new cap $E_p$ has a  supporting sphere centered at $p$, and
\[
E_p\subset B(e,2\delta).
\]

This reflection preserves the relevant unit directions. Indeed, if $z\in E_p'$ and $
q=e+p-z\in E_p$,
then
\begin{equation}\label{direction}
p-q=z-e.
\end{equation}
Thus the unit vector determined by the old pair $(e,z)$ is also determined by the new pair $(q,p)$. Moreover, all the points $p\in\mathcal {P}(E,\delta)$ form a $\delta$-net of $E$, so every direction represented by the old cap and its center is retained. Let $\mathcal A(E,\delta)$ be the family of caps $E_p$.

Choose the rapidly decreasing scales
\[
\delta_k=2^{-2^{k^2}},\qquad k\ge1.
\]
The rapid decay of $\delta_k$ ensures that later generations are less important at the scales used to
estimate the dimension.

Start with
\[
P_0=\{0\},
\qquad
A_0=\{S^{n-1}\}.
\]
Inductively, define
\[
P_k=
\bigcup_{E\in A_{k-1}}
\mathcal P(E,\delta_k),
\]
and
\[
A_k=\bigcup_{E\in A_{k-1}}\mathcal A(E,\delta_k).
\]

Every cap in $A_k$ lies on a supporting sphere centered at a point in $P_k$. It is also contained in the $2\delta_k$-neighbourhood of the center of its parent cap, which belongs to $P_{k-1}$. Hence
\begin{equation}\label{P conditon}
P_{k+1}\subset P_{k-1}^{(2\delta_k)}.
\end{equation}
This is the main geometric feature of the construction: each generation is clustered near the
generation two steps earlier.

Since $\sum_k\delta_k<\infty$, repeated application of \eqref{P conditon}, separately along the even and odd generations, shows that all the sets $P_k$ are in a fixed bounded region. Hence
\[
P=\bigcup_{k\ge0}P_k
\]
is bounded. Set
\[
K=\overline P.
\]
Then $K$ is compact.

We next verify that $S^{n-1}\subset K-K$. For each cap $E$, write $e(E)$ for the center of its supporting unit sphere. Then $E-e(E)$ is a set of unit directions on the cap. Define
\[
D_k=\bigcup_{E\in A_k}\bigl(E-e(E)\bigr).
\]
From \eqref{direction}, we know that 
\[
D_k\supset-D_{k-1}.
\]
Since $D_0=S^{n-1}$, every $D_k$ is contained in $S^{n-1}$, and $-S^{n-1}=S^{n-1}$,
we obtain
\[
D_k=S^{n-1},
\qquad\text{for every }k\ge0.
\]

Fix $u\in S^{n-1}$. For each $k\ge1$, choose a cap $E\in A_{k-1}$ and $z\in E$ such that
\[
z-e(E)=u.
\]
There exists $p\in\mathcal P(E,\delta_k)\subset P_k$ with
\[
|p-z|\le\delta_k.
\]
Since $e(E),p\in P$, we have 
\[
u_k=p-e(E)\in P-P\subset K-K.
\]
Moreover,
\[
|u_k-u|=|p-z|\le\delta_k\rightarrow0.
\]
Thus every unit vector is a limit of points in $K-K$. $K-K$ is compact and hence closed. Therefore
\[
S^{n-1}\subset K-K.
\]

Now we estimate the lower box-counting dimension. Let $M_k$ denote the number of caps in $A_k$. Since each point in $P_k$ produces a new cap, we have $\#P_k\leq M_k$.

The initial sphere satisfies
\[
M_1\lesssim \delta_1^{-(n-1)}.
\]
For $k\ge2$, every cap in $A_{k-1}$ has radius $2\delta_{k-1}$, so it contains at most 
\[
C\left(\frac{\delta_{k-1}}{\delta_k}\right)^{n-1}
\]
points separated at distance $\delta_k$. Therefore
\[
M_k\lesssim M_{k-1}
\left(\frac{\delta_{k-1}}{\delta_k}\right)^{n-1}.
\]
Then we have
\begin{equation}\label{count}
\#P_k\le M_k\lesssim\delta_k^{-(n-1)}.
\end{equation}

We estimate the covering number of 
\[
P=\bigcup_{j\geq0}P_j=\left(\bigcup_{j=0}^kP_j\right)\bigcup P_{k+1} \bigcup \left(\bigcup_{j\geq k+2}P_j\right)
\]
by cubes with side length $r_k = \delta_k^\alpha$,
where $\alpha > 1$ will be chosen below.

We cover the first $k$ generations by counting their points, and cover the generation $k+1$ by using the clustering near $P_{k-1}$.

First, we explain why the later generations are less important for estimating the dimension.

By the construction, we have the condition \eqref{P conditon}, which implies every point in $\bigcup_{j\geq k+2}P_j$ lies within distance
\[
\varepsilon_k
=2\sum_{i\ge k+1}\delta_i
\le4\delta_{k+1}
\]
of $\bigcup_{j=0}^{k+1}P_j$. Thus, 
\[
P\subset\left(\bigcup_{j=0}^{k+1}P_j\right)^{(\varepsilon_k)}
\]

The fast decay of $\delta_k$ implies, for any fixed $\alpha$, 
\[
\frac{\varepsilon_k}{r_k}\simeq\frac{\delta_{k+1}}{\delta_k^\alpha}\simeq\delta_k^{2^{2k+1}-\alpha}\longrightarrow0, \textrm{ as } k\rightarrow\infty.
\]
Therefore, for sufficiently large $k$, $\varepsilon_k\leq r_k$, and we have
\[
\bigcup_{j=0}^{k+1}P_j\subset P\subset \left(\bigcup_{j=0}^{k+1}P_j\right)^{(r_k)}.
\]
By enlarging each covering cube by a fixed factor, we have
\begin{equation}\label{tail contained}
N_{r_k}(P)\simeq N_{r_k}\left(\bigcup_{j=0}^{k+1}P_j\right).
\end{equation}

For the first $k$ generations, we first estimate $N_{r_k}\left(\bigcup_{j=1}^{k}P_j\right)$ directly by counting the points. Using \eqref{count} and the fast decay of $\delta_j$,
\[
N_{r_k}\left(\bigcup_{j=0}^{k}P_j\right)\le\sum_{j=0}^{k}\#P_j
\lesssim \sum_{j=1}^{k}\delta_j^{-(n-1)}
\lesssim \delta_k^{-(n-1)}.
\]

For generation $k+1$, recall that \eqref{P conditon} gives
\[
P_{k+1}\subset
\bigcup_{p\in P_{k-1}}B(p,2\delta_k).
\]
Each ball $B(p,2\delta_k)\subset\mathbb{R}^n$ can be covered by at most
\[
C_n\left(\frac{\delta_k}{r_k}\right)^n
\]
cubes of side length $r_k$. Consequently,
\[
N_{r_k}(P_{k+1})
\lesssim\left(\#P_{k-1}\right)
\left(\frac{\delta_k}{r_k}\right)^n\lesssim \delta_{k-1}^{-(n-1)}\delta_k^{n(1-\alpha)}.
\]
Choose 
\[
\alpha=\frac{2n-1}{n}.
\]
Then we have
\begin{equation}\label{sum}
    N_{r_k}(P)\simeq N_{r_k}\left(\bigcup_{j=1}^{k+1}P_j\right)\lesssim\delta_k^{-(n-1)}
+\delta_{k-1}^{-(n-1)}\delta_k^{n(1-\alpha)}\lesssim\delta_{k-1}^{-(n-1)}\delta_k^{-(n-1)}.
\end{equation}

Finally, by the definition of box counting dimension, we have
\begin{equation}
    \underline{\dim}_{B}(K)=\underline{\dim}_{B} (P)\leq
\liminf_{k\to\infty}
\frac{\log N_{r_k}(P)}{|\log r_k|}\leq\liminf_{k\to\infty}\frac{\log \delta_{k-1}^{-(n-1)}+\log \delta_k^{-(n-1)}}{\left|\log \delta_k^{\frac{2n-1}{n}}\right|}.
\end{equation}
Since $\delta_k=2^{-2^{k^2}}$, we have 
\[
\frac{\log \delta_{k-1}^{-(n-1)}}{\left|\log \delta_k^{\frac{2n-1}{n}}\right|}=\frac{(n-1)2^{(k-1)^2}}{\frac{2n-1}{n}2^{k^2}}\longrightarrow0, \textrm{ as }k\rightarrow\infty,
\]
and
\[
\frac{\log \delta_{k}^{-(n-1)}}{\left|\log \delta_k^{\frac{2n-1}{n}}\right|}=\frac{(n-1)2^{k^2}}{\frac{2n-1}{n}2^{k^2}}= \frac{n(n-1)}{2n-1}
\]
Therefore,
\[
\underline{\dim}_{B}(K)\leq\frac{n(n-1)}{2n-1}.
\]
This completes the construction and the proof.
\end{proof} 
In conclusion, $\frac{n(n-1)}{2n-1}$ is the sharp lower bound for the lower box counting dimension of compact sets $K\subset \mathbb{R}^n$ satisfying $S^{n-1}\subset K-K$. 
    \end{proof}
\begin{rem}
    The same construction also gives
\[
\dim_H (K)=0.
\] 

Recall that
\[
P=\bigcup_{j\ge0}P_j,
\qquad K=\overline P.
\]

Each $P_j$ is finite, so $P$ is countable and has Hausdorff dimension zero. It remains to show that 
 $\dim_H(K \setminus P)=0$.

Fix $k\ge2$. Repeated application of the clustering property
\[
P_{j+1}\subset P_{j-1}^{(2\delta_j)}
\]
shows that the generations $P_{k+1},P_{k+3},\ldots$ are close to $P_{k-1}$, that is,
\[
\bigcup_{\ell\ge0}P_{k+1+2\ell}
\subset P_{k-1}^{(a_k)},
\qquad
a_k=2\sum_{\ell\ge0}\delta_{k+2\ell}\le4\delta_k.
\]
Similarly, the generations $P_{k+2},P_{k+4},\ldots$ are close to $P_k$, that is
\[
\bigcup_{\ell\ge0}P_{k+2+2\ell}
\subset P_k^{(b_k)},
\qquad
b_k=2\sum_{\ell\ge0}\delta_{k+1+2\ell}\le4\delta_{k+1}.
\]
Therefore we have
\begin{equation}\label{set condition}
\bigcup_{j\geq k+1}P_j\subset P_{k-1}^{(4\delta_k)}\cup P_{k}^{(4\delta_{k+1})}.
\end{equation}

We claim that
\[
K\setminus P \subset \overline{\bigcup_{j\geq k+1}P_j}.
\]
Indeed, for any $x\in K\setminus P$ and $k\geq 2$, since $x\in\overline{P}$, there exists a sequence of points $x_m\in P$ converging to $x$. Since the set $\bigcup_{j=0}^kP_j$ is finite and $x\notin \bigcup_{j=0}^kP_j$, we have
\[
d_k:=\min_{y\in\bigcup_{j=0}^kP_j}|x-y|>0.
\]
Hence, for all sufficiently large $m$,
\[
|x_m-x|<d_k
\quad\Longrightarrow\quad
x_m\in P\setminus \bigcup_{j=0}^kP_j
=\bigcup_{j\ge k+1}P_j.
\]
It follows that
\[
K\setminus P
\subset\overline{\bigcup_{j\ge k+1}P_j}.
\]
Taking closures in \eqref{set condition}, we obtain
\[
K\setminus P\subset \left(\bigcup_{p\in P_{k-1}}B(p,4\delta_k)\right)\bigcup\left(\bigcup_{p\in P_{k}}B(p,4\delta_{k+1})\right).
\]
Recall $\#P_{k}\lesssim\delta_{k}^{-(n-1)}$ and the definition in \eqref{H}. For any $\alpha>0$,
\[
\mathcal{H}^\alpha_{8\delta_k}(K\setminus P)\lesssim \#P_{k-1}\delta_k^{\alpha}+\#P_k\delta_{k+1}^\alpha\lesssim \delta_{k-1}^{-(n-1)}\delta_k^\alpha+\delta_{k}^{-(n-1)}\delta_{k+1}^\alpha
\]
For any $\alpha>0$, since
\[
\lim_{k\to\infty}\delta_k^{-(n-1)}\delta_{k+1}^\alpha=\infty,
\]
we conclude
\[
\dim_{H}(K)=0.
\]

\end{rem}

\end{document}